\documentclass[12pt]{amsart}
\usepackage{amssymb,amsmath}
\usepackage{amsfonts}

\newtheorem{theorem}{Theorem}

\newtheorem{corollary}{Corollary}

\begin{document}
\title{Distributive lattices and cohomology.}

\author[Tomasz Maszczyk]{Tomasz Maszczyk\dag}
\address{Institute of Mathematics, Polish Academy of Sciences,
\newline Sniadeckich 8,\  00--956 Warszawa, Poland
\vspace{3mm} \hspace{5mm}
\newline  Institute of Mathematics,
University of Warsaw,
\newline Banacha 2,\ 02--097 Warszawa, Poland}
\email{maszczyk@mimuw.edu.pl}

\thanks{\dag The author was partially supported by KBN grants N201 1770 33 and 115/E-343/SPB/6.PR UE/DIE 50/2005-2008.}

\begin{abstract}
A resolution of the intersection of a finite number of subgroups
of an abelian group by means of their sums is constructed,
provided the lattice generated by these subgroups is distributive.
This is used for detecting singularities of modules over Dedekind
rings. A generalized Chinese remainder theorem is derived as a
consequence of the above resolution. The Gelfand-Naimark duality
between finite closed coverings of compact Hausdorff spaces and
the generalized Chinese remainder theorem is clarified.
\end{abstract}

\maketitle

\section{Introduction} The Gelfand-Naimark duality identifies lattices of closed subsets
in compact Hausdorff spaces with lattices opposite to surjective
systems of quotients of unital commutative C*-algebras. Therefore,
given a finite collection $I_{0},\ldots,I_{n}$ of closed *-ideals
in a C*-algebra $A=C(X)$ of continuous functions on a compact
Hausdorff space $X$, it identifies coequalizers  in the category
of compact Hausdorff spaces ($V(I)\subset X$ is the zero locus of
the ideal $I\subset A=C(X)$)
\begin{align}
\bigcup_{\alpha=0}^{n}V(I_{\alpha})\leftarrow\coprod_{\alpha=0}^{n}
V(I_{\alpha})\leftleftarrows\coprod_{\alpha,\beta=0}^{n}V(I_{\alpha})\cap
V(I_{\beta})\label{glu}\end{align} with equalizers  in the
category of unital commutative C*-algebras
\begin{align}
A/\bigcap_{\alpha=0}^{n}I_{\alpha}\rightarrow\prod_{\alpha=0}^{n}
A/I_{\alpha}\rightrightarrows\prod_{\alpha,\beta=0}^{n}A/I_{\alpha}+I_{\beta}.\label{eq}\end{align}
In particular, finite families of closed *-ideals intersecting to
zero correspond to finite families of closed subsets covering $X$.
In general, lattices of closed *-ideals in commutative unital
C*-algebras are always distributive, since they are isomorphic by
the Gelfand-Naimark duality to lattices opposite to sublattices of
subsets. Therefore one can think about finite families of closed
subsets in a compact Hausdorff space as of finite subsets in a
distributive lattice of ideals. By Hilbert's Nullstellensatz one
can replace a compact Hausdorff space and its closed subsets by an
affine algebraic set $X$ over an algebraically closed field and
its algebraic subsets on one hand, and closed ideals in a
C*-algebra by radical ideals in the algebra ${\mathcal O}[X]$ of
polynomial functions on $X$, on the other hand. One can take also
a finite set of monomial ideals in a ring of polynomials over a
field \cite{bru} as well as a finite set of congruences in the
ring of integers and the family of corresponding ideals. In all
above cases the fact that the diagram (\ref{eq}) is an equalizer
is a consequence of distributivity of a corresponding lattice of
ideals, and in view of the last example can be regarded as a
generalized Chinese remainder theorem. More examples can be
obtained from the fact that every algebra of finite representation
type has distributive lattice of ideals \cite{rie} and the
property of having distributive lattice of ideals is Morita
invariant and open under deformations of finite dimensional
algebras \cite{sch}.

The aim of the present paper is to show that the generalized
Chinese theorem is a consequence of vanishing of first cohomology
of a canonical complex associated with a finite number of members
$I_{0},\ldots,I_{n}$ of a distributive lattice $L$ of subgroups of
an abelian group $A$. The respective vanishing theorem (Theorem 1)
depends only on that lattice. Since it is independent of the
ambient abelian group $A$, Theorem 1 is prior to the generalized
Chinese remainder theorem. It is also more general, since it
claims vanishing of the whole higher cohomology. This can be used
for computing some higher Ext's detecting singularity of modules
over Dedekind rings (Corollary 1).

\section{Distributive lattices and homological algebra}
In this section we consider lattices of subgroups of a given
abelian group with the intersection and the sum as the meet and
the join operations, respectively. As for a general lattice the
distributivity condition can be written in two equivalent forms:

\begin{itemize}
\item The sum distributes over the intersection
\begin{align}
P_{0}\cap(P_{1}+ P_{2})=(P_{0}\cap P_{1})+ (P_{0}\cap P_{2}).
\end{align}

\item The intersection distributes over the sum
\begin{align}
I_{0}+I_{1}\cap I_{2}=(I_{0}+I_{1})\cap (I_{0}+I_{2}).
\end{align}
\end{itemize}
The aim of this section is to explain \emph{homological} nature of
the first and \emph{cohomological} nature of the second form of
distributivity. Homological characterization of distributivity was
used by Zharinov in \cite{zha} to generalize famous
edge-of-the-wedge theorem of Bogolyubov. In the present paper we
prove a cohomological characterization of distributivity and
derive from it a generalized Chinese remainder theorem. In this
section we show also that both characterizations have consequences
for arithmetic.

\subsection{Homology} Let $P_{0},\ldots,P_{n}$ be a finite family
of members of  some fixed lattice $L$ of subgroups in an abelian
group $A$. We define a group of $q$-chains ${\rm
C}_{q}(P_{0},\ldots,P_{n})$ as a quotient of the direct sum
\begin{align}\bigoplus_{0\leq
\alpha_{0},\ldots,\alpha_{q}\leq n}P_{\alpha_{0}}\cap\ldots\cap
P_{\alpha_{q}}
\end{align}
by a subgroup generated by elements
\begin{align}p_{\alpha_{0},\ldots,\alpha',\ldots,\alpha'',\ldots,\alpha_{q}}+p_{\alpha_{0},\ldots,\alpha'',\ldots,\alpha',\ldots,\alpha_{q}},\ \ \ \
p_{\alpha_{0},\ldots,\alpha,\ldots,\alpha,\ldots,\alpha_{q}},
\end{align}
and boundary operators (in terms of representatives of elements of
quotient groups)
\begin{align}
\partial: {\rm C}_{q}(P_{0},\ldots,P_{n})\rightarrow
{\rm C}_{q-1}(P_{0},\ldots,P_{n}),\nonumber\\
(\partial
p)_{\alpha_{0}\ldots\alpha_{q-1}}=\sum_{\alpha_{q}}p_{\alpha_{0}\ldots\alpha_{q-1}\alpha_{q}}.
\end{align}
By a standard argument from homological algebra $\partial\circ
\partial=0$. We denote by ${\rm H}_{\bullet}(P_{0},\ldots,P_{n})$
the homology of the complex $({\rm
C}_{\bullet}(P_{0},\ldots,P_{n}),
\partial)$.
\begin{theorem}{\cite{zha}}
1) ${\rm H}_{0}(P_{0},\ldots,P_{n})=P_{0}+\cdots + P_{n}$,

2) If the lattice $L$ is distributive then ${\rm
H}_{q}(P_{0},\ldots,P_{n})=0$ for $q>0$,

3) If ${\rm H}_{1}(P_{0},P_{1},P_{2})=0$ for all
$P_{0},P_{1},P_{2}\in L$ then the lattice $L$ is distributive.
\end{theorem}

\vspace{3mm} The following corollary provides a \emph{homological
characterization} of distributivity of a lattice $L$.
\begin{corollary} The following conditions are equivalent.

1) $L$ is distributive,

2) For all $P_{0},\ldots,P_{n}\in L$ the canonical morphisms of
complexes
\begin{align}
{\rm C}_{\bullet}(P_{0},\ldots,P_{n})\rightarrow P_{0}+\ldots+
P_{n}\label{hores}
\end{align}
are quasiisomorphisms.
\end{corollary}

In particular, since ${\rm C}_{q}(P_{0},\ldots,P_{n})$ can be
identified with the direct sum $\bigoplus_{0\leq
\alpha_{0}<\ldots<\alpha_{q}\leq n}P_{\alpha_{0}}\cap\ldots\cap
P_{\alpha_{q}}$, the above corollary provides a \emph{homological}
resolution of the sum of subgroups $P_{0}+\ldots+ P_{n}$ by means
of their intersections $P_{\alpha_{0}}\cap\ldots\cap
P_{\alpha_{q}}$, $0\leq \alpha_{0}<\ldots<\alpha_{q}\leq n$,
provided the lattice $L$ is distributive.

\vspace{3mm} In \cite{aliII} authors introduce a class of so
called G*GCD rings, defined as such for which gcd($P_{1},P_{2}$)
exists for all finitely generated projective ideals $P_{1},P_{2}$.
This class includes GGCD rings, semihereditary rings, f.f. rings
(and hence flat rings), von Neumann regular rings, arithmetical
rings, Pr\"ufer domains and GGCD domains. For every such a ring
the class of finitely generated projective ideals is closed under
intersection \cite{aliII}. Therefore, for an arithmetical ring $R$
every sum $P_{0}+\ldots+ P_{n}$ of finitely generated projective
ideals in $A$ admits a canonical resolution (\ref{hores}) by
finitely generated projective modules, which implies the following
corollary.
\begin{corollary} Let $P_{0},\ldots,P_{n}$ be finitely generated projective ideals in an arithmetical ring $R$. Then
\begin{align}
{\rm Ext}^{q}_{R}(P_{0}+\ldots+ P_{n}, -) & =  {\rm H}^{q}({\rm
Hom}_{R}({\rm C}_{\bullet}(P_{0},\ldots,P_{n}),-)),\\
{\rm Tor}_{q}^{R}(P_{0}+\ldots+ P_{n}, -) & =  {\rm H}_{q}({\rm
C}_{\bullet}(P_{0},\ldots,P_{n})\otimes_{R}-).
\end{align}
\end{corollary}

\subsection{Cohomology} Let $I_{0},\ldots,I_{n}$ be a finite family of members of  some
fixed lattice $L$ of subgroups in an abelian group $A$. We define
a group of $q$-cochains ${\rm C}^{q}(I_{0},\ldots,I_{n})$ as a
subgroup of the product
\begin{align}\prod_{0\leq
\alpha_{0},\ldots,\alpha_{q}\leq
n}I_{\alpha_{0}}+\ldots+I_{\alpha_{q}}\end{align} consisting of
sequences $(i_{\alpha_{0}\ldots\alpha_{q}})$ which are completely
alternating with respect to indices
$\alpha_{0},\ldots,\alpha_{q}$, i.e.
\begin{align}i_{\alpha_{0},\ldots,\alpha',\ldots,\alpha'',\ldots,\alpha_{q}}+i_{\alpha_{0},\ldots,\alpha'',\ldots,\alpha',\ldots,\alpha_{q}}=0,\ \ \ \
i_{\alpha_{0},\ldots,\alpha,\ldots,\alpha,\ldots,\alpha_{q}}=0,
\end{align}
and coboundary operators
\begin{align}
d: {\rm C}^{q}(I_{0},\ldots,I_{n})\rightarrow
{\rm C}^{q+1}(I_{0},\ldots,I_{n}),\nonumber\\
(di)_{\alpha_{0}\ldots\alpha_{q+1}}=\sum_{p=0}^{q+1}(-1)^{p}i_{\alpha_{0}\ldots\widehat{\alpha_{p}}\ldots\alpha_{q+1}}.
\end{align}
By a standard argument from homological algebra $d\circ d=0$. We
denote by ${\rm H}^{\bullet}(I_{0},\ldots,I_{n})$ the cohomology
of the complex $({\rm C}^{\bullet}(I_{0},\ldots,I_{n}), d)$.
\begin{theorem}
1) ${\rm H}^{0}(I_{0},\ldots,I_{n})=I_{0}\cap\ldots\cap I_{n}$,

2) If the lattice $L$ is distributive then ${\rm
H}^{q}(I_{0},\ldots,I_{n})=0$ for $q>0$,

3) If ${\rm H}^{1}(I_{0},I_{1},I_{2})=0$ for all
$I_{0},I_{1},I_{2}\in L$ then the lattice $L$ is distributive.
\end{theorem}

\emph{Proof.} Let us note first that ${\rm
C}^{q}(I_{0},\ldots,I_{n})$ can be identified with the product
$\prod_{0\leq \alpha_{0}<\ldots<\alpha_{q}\leq
n}I_{\alpha_{0}}+\cdots+ I_{\alpha_{q}}$.

1) Since the difference $i_{\beta}-i_{\alpha}$ is alternating with
respect to the indices $\alpha, \beta$ we have
\begin{align}
{\rm H}^{0}(I_{0},\ldots,I_{n})={\rm ker }(\prod_{0\leq\alpha\leq
n }I_{\alpha} & \rightarrow \prod_{0\leq\alpha<\beta\leq n
}I_{\alpha}+I_{\beta} ),\nonumber\\
(i_{\alpha}) & \mapsto (i_{\beta}-i_{\alpha})\nonumber
\end{align}
which consists of constant sequences $(i_{\alpha}=i\ \mid i\in
I_{0}\cap\ldots\cap I_{n})$.

2) For $q>0$ induction on $n$. For $n=0$ obvious. Inductive step:
Consider $(I_{0},\ldots,I_{n})$ for $n>0$. Then every $q$-cochain
$(i_{\alpha_{0}\ldots\alpha_{q}},
i_{\alpha_{0}\ldots\alpha_{q-1}n})$, for $q> 0$, can be identified
with a sequence consisting of elements
\begin{align}
i_{\alpha_{0}\ldots\alpha_{q}}  \in
I_{\alpha_{0}}+\ldots+I_{\alpha_{q}},\ \ {\rm for}\ \ 0\leq
\alpha_{0}<\ldots<\alpha_{q}\leq n-1,\\
i_{\alpha_{0}\ldots\alpha_{q-1}n}  \in
I_{\alpha_{0}}+\ldots+I_{\alpha_{q-1}}+I_{n},\ \ {\rm for}\ \
0\leq \alpha_{0}<\ldots<\alpha_{q-1}\leq n-1.\nonumber
\end{align}
This is a cocycle iff
\begin{align}
\sum_{p=0}^{q+1}(-1)^{p}i_{\alpha_{0}\ldots\widehat{\alpha_{p}}\ldots\alpha_{q+1}}=0,\label{1}\end{align}
for all $0\leq \alpha_{0}<\ldots<\alpha_{q+1}\leq n-1$, and
\begin{align}\sum_{p=0}^{q}(-1)^{p}i_{\alpha_{0}\ldots\widehat{\alpha_{p}}\ldots\alpha_{q}n}+(-1)^{q+1}i_{\alpha_{0}\ldots\alpha_{q}}=0.\label{2}
\end{align}
for all $0\leq \alpha_{0}<\ldots<\alpha_{q}\leq n-1$.

 By the inductive hypothesis ${\rm H}^{q}(I_{0},\ldots,I_{n-1})=0$
for $q>0$. Then (\ref{1}) implies that for all $0\leq
\alpha_{0}<\ldots<\alpha_{q-1}\leq n-1$ there exist
$i_{\alpha_{0}\ldots\alpha_{q-1}}\in
I_{\alpha_{0}}+\ldots+I_{\alpha_{q-1}}$, such that for all $0\leq
\alpha_{0}<\ldots<\alpha_{q}\leq n-1$
\begin{align}
i_{\alpha_{0}\ldots\alpha_{q}}=\sum_{p=0}^{q}(-1)^{p}i_{\alpha_{0}\ldots\widehat{\alpha_{p}}\ldots\alpha_{q}},\label{3}
\end{align}
hence by (\ref{2})
\begin{align}
\sum_{p=0}^{q}(-1)^{p}i_{\alpha_{0}\ldots\widehat{\alpha_{p}}\ldots\alpha_{q}n}+
\sum_{p=0}^{q}(-1)^{p+q+1}i_{\alpha_{0}\ldots\widehat{\alpha_{p}}\ldots\alpha_{q}}=0,
\end{align}
which can be rewritten as
\begin{align}
\sum_{p=0}^{q}(-1)^{p}(i_{\alpha_{0}\ldots\widehat{\alpha_{p}}\ldots\alpha_{q}n}-
(-1)^{q}i_{\alpha_{0}\ldots\widehat{\alpha_{p}}\ldots\alpha_{q}})=0.\label{4}
\end{align}
For $q=1$ we know by already proven point 2) of the theorem that
${\rm
H}^{q-1}(I_{0}+I_{n},\ldots,I_{n-1}+I_{n})=(I_{0}+I_{n})\cap\ldots\cap(I_{n-1}+I_{n})$
which is equal to $I_{0}\cap\ldots\cap I_{n-1}+I_{n}$, since the
lattice $L$ is distributive. Therefore by (\ref{4}) for $q=1$
there exist $i\in I_{0}\cap\ldots\cap I_{n-1}$ and $i_{n}\in
I_{n}$ such that for all $0\leq \alpha\leq n-1$
\begin{align}
i_{\alpha n}+i_{\alpha}=i+i_{n},
\end{align}
hence
\begin{align}
i_{\alpha n}=i_{n}-(i_{\alpha}-i).\label{6}
\end{align}
Equations (\ref{3})  for $q=1$, which reads as
\begin{align}
i_{\alpha
\beta}=i_{\beta}-i_{\alpha}=(i_{\beta}-i)-(i_{\alpha}-i),
\end{align}
and (\ref{6}) together mean that $i=(i_{\alpha\beta}, i_{\alpha
n}\ \mid 0\leq \alpha<\beta\leq n-1)\in {\rm C}^{1}(I_{0},\ldots,
I_{n})$ is coboundary of $(i_{\alpha}-i, i_{n}\ \mid 0\leq
\alpha\leq n-1)\in {\rm C}^{0}(I_{0},\ldots, I_{n})$ which proves
that ${\rm H}^{1}(I_{0},\ldots,I_{n})=0$.

For $q>1$ by the inductive hypothesis ${\rm
H}^{q-1}(I_{0}+I_{n},\ldots,I_{n-1}+I_{n})=0$, hence (\ref{4})
implies that for all $0\leq \alpha_{0}<\ldots<\alpha_{q-2}\leq
n-1$ there exist $i_{\alpha_{0}\ldots\alpha_{q-2}n}\in
(I_{\alpha_{0}}+I_{n})+\ldots+(I_{\alpha_{q-2}}+I_{n})=I_{\alpha_{0}}+\ldots+I_{\alpha_{q-2}}+I_{n}$,
such that for all $0\leq \alpha_{0}<\ldots<\alpha_{q-1}\leq n-1$
\begin{align}
i_{\alpha_{0}\ldots\alpha_{q-1}n}-
(-1)^{q}i_{\alpha_{0}\ldots\alpha_{q-1}}=\sum_{p=0}^{q-1}(-1)^{p}i_{\alpha_{0}\ldots\widehat{\alpha_{p}}\ldots\alpha_{q-1}n},
\end{align}
which can be rewritten as
\begin{align}
i_{\alpha_{0}\ldots\alpha_{q-1}n}=\sum_{p=0}^{q-1}(-1)^{p}i_{\alpha_{0}\ldots\widehat{\alpha_{p}}\ldots\alpha_{q-1}n}+
(-1)^{q}i_{\alpha_{0}\ldots\alpha_{q-1}}.\label{5}
\end{align}
Equations (\ref{3}) and (\ref{5}) together mean that
$(i_{\alpha_{0}\ldots\alpha_{q}},
i_{\alpha_{0}\ldots\alpha_{q-1}n})$ is coboundary of
$(i_{\alpha_{0}\ldots\alpha_{q-1}},
i_{\alpha_{0}\ldots\alpha_{q-2}n})$, hence ${\rm
H}^{q}(I_{0},\ldots,I_{n})=0$ for $q>1$.

3) We have to prove that for all $I_{0},I_{1},I_{2}\in L$
\begin{align}
(I_{0}+I_{1})\cap (I_{0}+I_{2})=I_{0}+I_{1}\cap I_{2}.
\end{align}
The inclusion $(I_{0}+I_{1})\cap (I_{0}+I_{2})\supset
I_{0}+I_{1}\cap I_{2}$ is obvious. To prove the opposite inclusion
take $i\in (I_{0}+I_{1})\cap (I_{0}+I_{2})$. It can be written in
two ways as
\begin{align}
i=i_{01}+i_{12}',\ \ \ {\rm where}\ \ \ i_{01}\in I_{0}\subset
I_{0}+I_{1},\ \ i_{12}'\in I_{1}\subset I_{1}+I_{2},\label{01}
\end{align}
\begin{align}
i=i_{02}+i_{12}'',\ \ \ {\rm where}\ \ \ i_{02}\in I_{0}\subset
I_{0}+I_{2},\ \ i_{12}''\in I_{2}\subset I_{1}+I_{2}.\label{02}
\end{align}
Define $i_{12}:=i_{12}'-i_{12}''$. Subtracting (\ref{02}) from
(\ref{01}) we get the cocycle condition
\begin{align}
i_{12}-i_{02}+ i_{01}=0.\label{coc}
\end{align}
Since ${\rm H}^{1}(I_{0},I_{1},I_{2})=0$ (\ref{coc}) implies that
there exist $i_{\alpha}\in I_{\alpha}$, $\alpha=0,1,2$, such that
\begin{align}
i_{\alpha\beta}=i_{\beta}- i_{\alpha},\label{cob}
\end{align}
in particular
\begin{align}
i_{12}'-i_{12}''=i_{12}=i_{2}-i_{1},
\end{align}
hence
\begin{align}
i_{1}+ i_{12}'=i_{2}+ i_{12}''.\label{'}
\end{align}
Since $i_{0}\in I_{0}$ and by (\ref{01}) $i_{01}\in I_{0}$
(\ref{cob}) implies that
\begin{align}i_{1}=i_{0}+i_{01}\in I_{0}.\label{I0}\end{align}

By (\ref{01}) (resp. (\ref{02})) the left (resp. right) hand side
of (\ref{'}) belongs to $I_{1}$ (resp. $I_{2}$), hence
\begin{align}
i_{1}+ i_{12}'\in I_{1}\cap I_{2}.\label{I12}
\end{align}
Finally, by (\ref{01}), (\ref{I0}) and (\ref{I12})
\begin{align}
i=i_{01}+ i_{12}'=(i_{01}-i_{1})+ (i_{1}+i_{12}')\in
I_{0}+I_{1}\cap I_{2}. \ \Box
\end{align}

\vspace{3mm} The following corollary provides a
\emph{cohomological characterization} of distributivity of a
lattice $L$.
\begin{corollary} The following conditions are equivalent.

1) $L$ is distributive,

2) For all $I_{0},\ldots,I_{n}\in L$ the canonical morphisms of
complexes
\begin{align}
I_{0}\cap\ldots\cap I_{n}\rightarrow {\rm
C}^{\bullet}(I_{0},\ldots,I_{n}),
\end{align}
are quasiisomorphisms.
\end{corollary}

In particular, since ${\rm C}^{q}(I_{0},\ldots,I_{n})$ can be
identified with the product $\prod_{0\leq
\alpha_{0}<\ldots<\alpha_{q}\leq n}I_{\alpha_{0}}+\cdots+
I_{\alpha_{q}}$,  the above corollary provides a
\emph{cohomological} resolution of the intersection of subgroups
$I_{0}\cap\ldots\cap I_{n}$ by means of their sums
$I_{\alpha_{0}}+\cdots+I_{\alpha_{q}}$, $0\leq
\alpha_{0}<\ldots<\alpha_{q}\leq n$, provided the lattice $L$ is
distributive.

\vspace{3mm}This fact together with the fact that each finitely
generated module over a Dedekind ring is a direct sum of
distributive modules (i.e. modules whose lattice of submodules is
distributive) \cite{tug} can be used for detecting singularities
of modules over Dedekind rings. First of all, in a non-singular
left module (i.e. left module without nonzero elements annihilated
by all essential left ideals) the intersection of injective
submodules is again injective \cite{wis}. Therefore, given
injective submodules $I_{0},\ldots,I_{n}$ in a left module $A$
over a ring $R$, the functors ${\rm Ext}^{q}_{R}(-,
I_{0}\cap\ldots\cap I_{n})$ for $q>0$ detect singularity of $A$.
These functors can be computed with use of our resolution whenever
every sum of injective submodules of a left $R$-module $A$ is
injective. The latter property characterizes left Noetherian left
hereditary rings \cite{he}, hence it holds for Dedekind rings.
Therefore we can apply Theorem 2 to obtain the following
corollary.

\begin{corollary} Let $I_{0},\ldots,I_{n}$ be injective submodules
in a distributive left $R$-module $A$ over a left Noetherian and
left hereditary ring $R$. Then
\begin{align}
{\rm Ext}^{q}_{R}(-, I_{0}\cap\ldots\cap I_{n})={\rm H}^{q}({\rm
Hom}_{R}(-,{\rm C}^{\bullet}(I_{0},\ldots,I_{n}))).
\end{align}
Therefore, if $A$ is a finitely generated and non-singular module
over a Dedekind ring $R$
\begin{align}{\rm H}^{q}({\rm Hom}_{R}(-,{\rm
C}^{\bullet}(I_{0},\ldots,I_{n})))=0\end{align} for $q>0$.
\end{corollary}

\section{Generalized Chinese
Remainder Theorem} As next application we will prove the following
generalized Chinese remainder theorem.
\begin{corollary}
For any finite family $ I_{0},\ldots,I_{n}$ of members of some
fixed distributive lattice $L$ of subgroups in an abelian group
$A$ the canonical diagram
\begin{align}
A/\bigcap_{\alpha=0}^{n}I_{\alpha}\rightarrow\prod_{\alpha=0}^{n}A/I_{\alpha}\rightrightarrows\prod_{\alpha,\beta=0}^{n}A/I_{\alpha}+I_{\beta}.\label{7}\end{align}
is an equalizer.
\end{corollary}

\emph{Proof.} Injectivity of the first arrow is obvious. Exactness of (\ref{7}) in the middle term is equivalent to
exactness in the middle term of the canonical complex
\begin{align}
A\stackrel{\pi}{\rightarrow}\prod_{0\leq\alpha\leq n
}A/I_{\alpha}\stackrel{\delta}{\rightarrow}\prod_{0\leq\alpha<\beta\leq
n}A/I_{\alpha}+I_{\beta},\label{8}\end{align} where
$\pi(a)=(a+I_{\alpha} \mid 0\leq\alpha\leq n )$,
$\delta(a_{\alpha}+I_{\alpha} \mid 0\leq\alpha\leq
n)=(a_{\beta}-a_{\alpha}+I_{\alpha}+I_{\beta}\ \mid
0\leq\alpha<\beta\leq n)$. We have
\begin{align}
{\rm ker}\ \delta=(a_{\alpha}+I_{\alpha}\mid
a_{\beta}-a_{\alpha}\in I_{\alpha}+I_{\beta}),
\end{align}
hence $(i_{\alpha\beta}:=a_{\beta}-a_{\alpha}\mid
0\leq\alpha<\beta\leq n)$ is a cocycle in ${\rm
C}^{1}(I_{0},\ldots,I_{n})$. By Theorem 1 there exist
$i_{\alpha}\in I_{\alpha}$ such that
$a_{\beta}-a_{\alpha}=i_{\beta}-i_{\alpha}$. Let
$a:=a_{\alpha}-i_{\alpha}=a_{\beta}-i_{\beta}$. Then
$(a_{\alpha}+I_{\alpha}\mid 0\leq\alpha\leq n)=\pi(a)$, which
proves exactness of (\ref{8}) in the middle term. $\Box$

\vspace{3mm} \textbf{Remark.} It is well known that if all
$I_{\alpha}$'s are pairwise coprime ideals in a unital associative
ring $A$, i.e. $I_{\alpha}+I_{\beta}=A$ for $\alpha\neq\beta$,
then the diagram (\ref{7}) is an equalizer and
$I_{1}\cap\ldots\cap I_{n}=\sum_{\sigma}I_{\sigma(1)}\ldots
I_{\sigma(n)}$, where $\sigma$'s are sufficiently many
permutations of $\{ 1,\ldots, n\}$. These facts are independent of
distributivity of the lattice of ideals. Therefore Corollary 5
(essentialy present already in \cite{mac}, next rediscovered and
generalized many times, e.g. \cite{bal, ach, vag, bru}) should be
understood as a generalization of the Chinese remainder theorem to
the non-coprime case, for which distributivity of the lattice of
ideals is a sufficient condition. In fact, the lattice of left
ideals in a (unital associative) ring is distributive iff the
above generalized Chinese remainder theorem holds for such ideals
\cite{ach}. Therefore in the commutative case there is \emph{``one
necessary and sufficient condition that places the theorem in
proper perspective. It states that the Chinese remainder theorem
holds in a commutative ring if and only if the lattice of ideals
of the ring is distributive''} \cite{rot}. The aim of this section
was to show how lattice theory communicates with modular
arithmetic through homology theory.

\section{Noncommutative Topology}
Finite families of closed subsets covering a topological space are
important for the Mayer-Vietoris principle in sheaf cohomology
with supports and topological K-theory. Since by the
Gelfand-Naimark duality gluing of a compact Hausdorff space $X$
from finite number of compact Hausdorff pieces is equivalent to a
generalized Chinese remainder theorem (\ref{eq}) for closed
*-ideals in a commutative unital C*-algebra $C(X)$, one is tempted
to define a \emph{``noncommutative closed covering of a noncommutative
space dual to an associative C*-algebra $A$''} as a finite
collection of closed *-ideals intersecting to zero and generating
a distributive lattice \cite{cal}, \cite{haj}.

In \cite{haj}
authors focus on the combinatorial side of such gluing in terms of the poset structure on $X$ induced by such a covering. This poset structure has its own topology (Alexandrov topology), drastically different from the original compact Hausdorff one.  After fixing an order of the finite
closed covering, they represent the distributive lattice generated
by these originally closed (now Alexandrov open) subsets as a homomorphic image of the
free distributive lattice generated by the same finite set of
generators.  Next, authors pull-back quotient C*-algebras $A/I$ to that free
lattice and view the resulting surjective system of
quotient algebras as \emph{flabby sheaf} of C*-algebras on the
Alexandrov topology corresponding to that free lattice. Finally,
they formulate the Gelfand-Naimark duality between ordered
coverings of compact Hausdorff spaces by $N$ closed sets and
flabby sheaves of commutative unital C*-algebras on the Alexandrov
topology corresponding to the free distributive lattice with $N$
generators.

The aim of the present chapter is to avoid the auxiliary Alexandrov topology.
In fact, creating new topology by declaring old
closed subsets to be new opens is not necessary. The reason is that
there is no need to see the generalized Chinese remainder theorem
as the sheaf condition. The following definitions introduce a notion, which replaces sheaves when finite closed coverings replace open coverings.

\vspace{3mm} \paragraph{\textbf{Definition.}}    For any  topological space $X$ we
define a category of functors ${\mathcal P}$ (we call them
\emph{patterns}) from the lattice of \emph{closed subsets} of $X$
to the category of sets (abelian groups, rings, algebras etc)
satisfying the following \emph{unique gluing property} with
respect to finite closed coverings $(C_{0},\ldots,C_{n})$ of
closed subsets $C=C_{0}\cup\ldots\cup C_{n}\subset X$, which
demands that all canonical diagrams
\begin{align}
{\mathcal P}(C)\rightarrow\prod_{\alpha=0}^{n} {\mathcal
P}(C_{\alpha})\rightrightarrows\prod_{\alpha,\beta=0}^{n}{\mathcal
P}(C_{\alpha}\cap C_{\beta})
\end{align}
are equalizers.

\vspace{3mm} \paragraph{\textbf{Definition.}} We call a pattern ${\mathcal P}$ on $X$ \emph{global} if for any
closed subset $C\subset X$ the restriction morphism ${\mathcal
P}(X)\rightarrow {\mathcal P}(C)$ is surjective.

\vspace{3mm} \paragraph{\textbf{Definition.}} For a continuous map $f:X\rightarrow Y$ the preimage $f^{-1}(D)$
of any closed subset $D\subset Y$ is closed in $X$ and $f$ defines
the \emph{direct image} functor $f_{\ast}$ of patterns:
$(f_{\ast}{\mathcal P})(D):={\mathcal P}(f^{-1}(D))$.

We call (\emph{globally}) \emph{algebraized space} a pair
consisting of a topological space $X$ equipped with a (global)
pattern ${\mathcal A}_{X}$ of algebras.

\vspace{3mm} \paragraph{\textbf{Definition.}} A morphism of (globally) algebraized spaces consists of a
continuous map of topological spaces $f:X\rightarrow Y$ and a
morphism of patterns of algebras ${\mathcal A}_{Y}\rightarrow
f_{\ast}{\mathcal A}_{X}$.

\vspace{3mm}
In this framework the aforementioned Gelfand-Naimark
duality between gluing  of compact Hausdorff spaces and the
generalized Chinese remainder theorem for C*-algebras reads now
as follows.
\begin{theorem}
The Gelfand-Naimark duality induces a full embedding of the
category opposite to unital commutative C*-algebras equipped with
lattices of closed *-ideals into the category of compact
Hausdorff globally algebraized spaces.
\end{theorem}

Note that in the above theorem the Gelfand-Naimark duality between
pairs $(A, L)$ and $(X, {\mathcal A}_{X})$  dualizes a unital
commutative  C*-algebra $A$ to a compact Hausdorff space $X$ and
the lattice $L$ of closed *-ideals in $A$ to a global pattern of
algebras ${\mathcal A}_{X}$.

Note that every lattice of closed *-ideals in a C*-algebra is
distributive. Therefore, according to the general ideology of
noncommutative topology,  a pair consisting of a unital
associative C*-algebra $A$ and a lattice $L$ of closed *-ideals in
$A$ should be regarded as a \emph{``noncommutative compact Hausdorff
globally algebraized space''}.

\subsection{C*-algebras and patterns}
\subsubsection{Continuous fields of C*-algebras}  In functional analysis of function C*-algebras, in
opposite to algebraic geometry, the notion of sheaf plays no a
significant role. The appropriate replacement is then the notion
of sections of \emph{continuous fields of C*-algebras} \cite{dix}.
They are contravariant functors on the category of closed subsets
transforming closed embeddings into surjective restriction
homomorphisms. It is easy to observe that they satisfy the unique
gluing property with respect to finite closed coverings of closed
subsets , i.e. they are patterns in our terminology. This property
was used in computation of $K$-theory of an important class of
Toeplitz algebras on Lie groups, with use of the Mayer-Vietoris
sequence \cite{par}.

\subsubsection{Continuous functions vanishing at infinity} Another
example of patterns arising in theory of function C*-algebras
comes from continuous functions vanishing at infinity. For any
locally compact space one has a non-unital C*-algebra
$\mathcal{C}_{0}(X)$ of continuous functions vanishing at
infinity. It is widely accepted that $\mathcal{C}_{0}(X)$ is an
appropriate C*-algebraic replacement of the locally compact space
$X$, mostly in view of the Gelfand-Naimark duality in the
unital-versus-compact case. A beautiful part of functional
analysis was created to extend the Gelfand-Naimark duality in this
way. However, an idea that relaxing compactness to local compactness can be dualized to forgetting about the unit of the C*-algebra is not true if one wants to preserve the usual relation between continuous functions and topology. 

First, about C*-algebras and locality. Although continuous
functions form a sheaf under restriction to open subsets, the
vanishing at infinity property does not survive the restriction.
This means that given two open subsets $U\subset V\subset X$ there
is no a restriction homomorphism from $\mathcal{C}_{0}(V)$ to $\mathcal{C}_{0}(U)$. Strange enough, there is a well defined injective
homomorphism of non-unital algebras in the opposite direction
$\mathcal{C}_{0}(U)\rightarrow \mathcal{C}_{0}(V)$, given by the extension
by zero. Moreover, given open subsets $U_{0},\ldots, U_{n}$ one
has a canonical equalizer diagram
\begin{align}
\mathcal{C}_{0}(U)\rightarrow\prod_{\alpha=0}^{n}
\mathcal{C}_{0}(U_{\alpha})\rightrightarrows\prod_{\alpha,\beta=0}^{n}\mathcal{C}_{0}(U_{\alpha}\cup
U_{\beta} ),\label{zeroexteq}\end{align} whose arrows are defined
as collections of extensions by zero with respect to inclusions
$U=U_{0}\cap\ldots\cap U_{n}\subset U_{\alpha}$,
$U_{\alpha}\subset U_{\alpha}\cup U_{\beta}$ and $U_{\beta}\subset
U_{\alpha}\cup U_{\beta}$.

 The \v{C}ech-Stone compactification $X\hookrightarrow
\beta X$ and the Gelfand-Naimark duality put the problem of geometric description of the extension by zero into the right perspective. The extension by zero $\mathcal{C}_{0}(U)\rightarrow
\mathcal{C}_{0}(V)$ is equivalent to a surjective restriction
homomorphism of unital quotient algebras
\begin{align}
\mathcal{C}(\beta X \setminus U)=\mathcal{C}(\beta
X)/\mathcal{C}_{0}(U)\rightarrow \mathcal{C}(\beta
X)/\mathcal{C}_{0}(V)=\mathcal{C}(\beta X \setminus V),
\end{align}
when we regard $\mathcal{C}_{0}(U)$ as a closed  *-ideal in the
C*-algebra $\mathcal{C}(\beta X)\cong\mathcal{C}_{b}(X)$ of
continuous functions on $\beta X$ isomorphic to the C*-algebra of
bounded functions on $X$. This restriction homomorphism is
Gelfand-Naimark dual to the closed inclusion $C:=\beta X \setminus
U\subset \beta X \setminus V=: D$, and makes (\ref{zeroexteq}) the
equalizer diagram verifying  the pattern property of the assignment
$C\mapsto \mathcal{I}(C):=\mathcal{C}_{0}(\beta X
\setminus C)$ on the finite closed covering
\begin{align}\beta X\setminus U = \bigcup_{\alpha=0}^{n}(\beta X\setminus U_\alpha).\end{align}

This means that the pattern $\mathcal{I}$ is a (sub)pattern of
ideals in the constant pattern $\mathcal{C}(\beta X)$ of algebras.
The pattern $C\mapsto \mathcal{C}(C)$ is then the
pattern of quotient algebras.

\subsubsection{Pattern cohomology on topological spaces} Patterns admit an analog of the \v{C}ech cohomology with respect to finite closed coverings. Assume that there is given such a covering $X=C_{0}\cup\ldots\cup C_{n}$ of a topological space $X$ and a pattern ${\mathcal P}$. Mimicking the \v{C}ech complex construction we define the \emph{pattern cohomology}
  
\begin{align}
{\rm H}^{p}(C_{0},\ldots, C_{n}; {\mathcal P}):={\rm H}^{p}(\prod_{i_{0} <\ldots < i_{\bullet}}{\mathcal P}(C_{i_{0}}\cap\ldots\cap C_{i_{\bullet}})).
\end{align}
The pattern property implies that
 \begin{align}
{\rm H}^{0}(C_{0},\ldots, C_{n}; {\mathcal P})={\mathcal P}(X).
\end{align}
If the pattern takes values in a distributive lattice of subgroups of an abelian group $A$ in such a way that
\begin{align}
{\mathcal P}(C_{i_{0}}\cap\ldots\cap C_{i_{p}}) & ={\mathcal P}(C_{i_{0}})+\ldots + {\mathcal P}(C_{i_{p}})\label{intersect}\\
{\mathcal P}(C_{i_{0}}\cup\ldots\cup C_{i_{p}}) & ={\mathcal P}(C_{i_{0}})\cap\ldots \cap {\mathcal P}(C_{i_{p}}).\label{union}
\end{align}
then by Theorem ? we obtain for $p>0$

\begin{align}
{\rm H}^{p}(C_{0},\ldots, C_{n}; {\mathcal P})=0.
\end{align}
In particular, for the constant pattern $A(C):=A$
\begin{align}
{\rm H}^{p}(C_{0},\ldots, C_{n}; A)=\left\{ \begin{array}{ll}
                                                   A &  {\rm if}\ p=0,\\
                                                   0 &  {\rm if}\ p\neq 0
                                                   \end{array}\right.\end{align}
and for its subpattern ${\mathcal I}$  taking values in a distributive lattice of subgroups of $A$, satisfying (\ref{intersect}) and (\ref{union}), and such that ${\mathcal I}(X)=0$
\begin{align}
{\rm H}^{p}(C_{0},\ldots, C_{n}; {\mathcal I})=0\end{align}
for all $p$. The short exact sequence of patterns
\begin{align}
0\rightarrow {\mathcal I}\rightarrow A\rightarrow A/{\mathcal I}\rightarrow 0\end{align}
induces then a long exact sequence of pattern cohomology, which implies that
\begin{align}
{\rm H}^{p}(C_{0},\ldots, C_{n}; A/{\mathcal I})=\left\{ \begin{array}{ll}
                                                   A &  {\rm if}\ p=0,\\
                                                   0 &  {\rm if}\ p\neq 0
                                                   \end{array}\right.\end{align}
This implies that the cohomological behavior of taking remainders modulo ideals (restrictions to closed subsets) of an arithmetical ring $A$ expressed in terms of the globally algebraized space structure defined on the Zariski topology of ${\rm Spec}(A)$ is similar to the cohomological behavior of localizations (restrictions to open subsets) of $A$ expressed in terms of the locally ringed space structure on ${\rm Spec}(A)$.
\subsubsection{Sheaves versus patterns}
 Due to Cartan, Leray's \emph{``faisceaux continus"} on locally compact spaces are equivalent to sheaves. Actually, given a sheaf $\mathcal{F}$ on a locally compact space $X$ one can assign to every closed subset $C\subset X$ the stalk of $\mathcal{F}$ along $C$. This assignment is different from our ``pattern". For a sheaf of continuous functions the stalk at a point consists of germs, while the evaluation of the pattern on a point consists of values. If the space $X$ is not discrete the kernel of the surjective evaluation map from the stalk to the set of values is usually big.

\end{document}